\documentstyle[12pt]{article}

\newtheorem{lem}{Lemma}
\newtheorem{thm}{Theorem}

\newtheorem{propref}{Proposition}
\newtheorem{cor}{Corollary}

\newcommand{\be}{\begin{equation}}
\newcommand{\ee}{\end{equation}}

\def\real{\hbox{\rm\setbox1=\hbox{I}\copy1\kern-.45\wd1 R}}

\def\R{\hbox{\rm\setbox1=\hbox{I}\copy1\kern-.45\wd1 R}}

\def\prob{\hbox{\rm\setbox1=\hbox{I}\copy1\kern-.45\wd1 P}}

\def\pr{\scriptsize{\hbox{\rm\setbox1=\hbox{I}\copy1\kern-.45\wd1 P}}}

\def\r{\scriptsize{\hbox{\rm\setbox1=\hbox{I}\copy1\kern-.45\wd1 R}}}

\def\Zint{{\mathchoice{\setbox1=\hbox{\sf Z}\copy1\kern-.75\wd1\box1}
         {\setbox1=\hbox{\sf Z}\copy1\kern-.75\wd1\box1}
         {\setbox1=\hbox{\footnotesize\sf Z}\copy1\kern-.75\wd1\box1}
         {\setbox1=\hbox{\footnotesize\sf Z}\copy1\kern-.75\wd1\box1}}}

\def\pint{\hbox{\rm\setbox1=\hbox{I}\copy1\kern-.55\wd1 N}}

\def\bx{{\bf x}}

\def\by{{\bf y}}
\def\bv{{\bf v}}

\def\bw{{\bf w}}

\newcommand{\deq}{\stackrel{\scriptscriptstyle\triangle}{=}}

\def\len{\mathop{\rm len}}

\def\clos{\mathop{\rm clos}}

\newcommand{\F}{\mbox{$\cal F$}}

\newcommand{\E}{\mbox{$\cal E$}}

\newcommand{\mh}{\mbox{$\hat{m}$}}

\newcommand{\hmu}{\hat{\mu}}
\newcommand{\hnu}{\hat{\nu}}

\newcommand{\tm}{\tilde{m}}

\def\B{{\cal B}}

\def\A{{\cal A}}

\def\H{{\cal H}}

\def\L{\mathop{\rm \Lambda}}

\begin{document}

\bibliographystyle{plain}

\thispagestyle{empty}
\setcounter{page}{0}

\vspace {2cm}

{\Large Guszt\'{a}v Morvai, Sanjeev R. Kulkarni, Andrew B. Nobel:}

 \vspace {1cm}

 \centerline {\Huge   Regression Estimation}

\vspace {1cm}

 \centerline {\Huge  from an
       Individual Stable Sequence.}

\vspace {1cm}

{\Large Appeared in :  Statistics Vol. 33, pp. 99-118, 1999.}

\vspace {2cm}

\begin{abstract}

We consider univariate regression estimation
from an individual (non-random) sequence 
$(x_1,y_1),(x_2,y_2), \ldots \in \real \times \real$,
which is stable in the sense that
for each interval $A \subseteq \real$,
(i) the limiting relative frequency of $A$ under
$x_1, x_2, \ldots$ is governed by an unknown probability
distribution $\mu$, and (ii) the limiting average of those $y_i$ with
$x_i \in A$ is governed by an unknown regression function $m(\cdot)$.

A computationally simple scheme for estimating $m(\cdot)$
is exhibited, and is shown to be $L_2$ consistent for 
stable sequences $\{(x_i,y_i)\}$ such that $\{y_i\}$ is bounded
and there is a known upper bound for the variation of 
$m(\cdot)$ on intervals of the form $(-i,i]$, $i \geq 1$.
Complementing this positive result, it is shown that 
there is no consistent estimation scheme for the family
of stable sequences whose regression functions
have finite variation, even under the restriction that
$x_i \in [0,1]$ and $y_i$ is binary-valued.

\end{abstract}


\newpage

\section{Introduction}

Individual numerical sequences (binary and real-valued) have played 
an important role in the theory of data compression and computational
complexity.
The theory of lossless data compression developed
by Ziv and Lempel \cite{LempelZiv76}, Ziv \cite{Ziv78}, and the 
complexity theory of Kolmogorov \cite{Kolmogorov65, Kolmogorov68} and
Chaitin \cite{Chaitin66} 
are both formulated within a purely deterministic framework 
that is built around individual sequences.
Subsequent work in these areas has considered
useful notions of randomness, compressibility, and predictability.
More recently, individual sequences have been studied in the context of
statistical learning theory.  
In spite of the above research, there has been little consideration of
individual sequences in the context of classical statistical
estimation.

It is common in statistics to treat data, for the purposes of
analysis, as a sequence of (typically independent) identically
distributed random variables.  
This stochastic point of view collapses when one is faced with 
a {\em particular} collection of data, which is 
a fixed sequence of numbers or vectors
from which we hope to learn something about the state of nature.

It is natural then to (re)formulate some classical statistical
problems in terms of individual sequences.  
We concern ourselves here with the important problem of
regression estimation.
In the common statistical setting one is given $n$ independent
replicates $(X_1,Y_1), \ldots, (X_n,Y_n)$ of a jointly distributed
pair $(X,Y) \in \real \times \real$, and asked to find an estimate
of the regression function $m(x) = E[Y|X=x]$.
Justification for estimation of $m(x)$ comes from the fact that
it minimizes $E(h(X) - Y)^2$ over all functions $h(\cdot)$ of $X$.
Thus $m(\cdot)$ is the least squares estimate of $Y$ given $X$.

In this paper we present and analyze a simple regression estimation
procedure that is applicable in a purely deterministic setting.  
By applying our estimates to individual sample paths, 
we easily establish their almost-sure consistency for 
ergodic processes having suitable one-dimensional
distributions (the dependence structure of the process is unimportant).
The approach and results of this paper are motivated by, 
and closely related to, recent results of \cite{NobMorKul96} on density
estimation from individual sequences.

For independent and weakly dependent stochastic data, 
a variety of estimation schemes have been proposed,
including procedures based on histograms, 
kernels, neural networks, orthogonal series, wavelets,
and nearest neighbors.
For a description of some of these methods  see, for example,
Gy\"orfi, H\"ardle, Sarda, and Vieu \cite{GHSV89}, Roussas \cite{Rou91},
Devroye, Krzyzak \cite{DevKrzy89} and the references therein.
Kulkarni and Posner \cite{KulPos95} 
studied nearest neighbor regression estimates
in the case where $x_1, x_2, \ldots$ are deterministic, but
$Y_1, Y_2, \ldots$ are random and conditionally independent
given the $x_i$'s.
Yakowitz {\em et al.} \cite{YGKM97} considered a family of 
truncated histogram regression 
estimates for processes with vector-valued covariates.
For each constant $L > 0$ they exhibit a sequence of
estimates that is almost surely pointwise consistent
for every ergodic process $\{(X_i,Y_i)\}$ whose regression function
satisfies a Lipschitz condition of the form
$|m(x) - m(y)| \leq L ||x - y||$.
In practice, the constant $L$ is known and fixed in advance of the data.
Related work has been done in the
area of nonparametric forecasting for a stationary process $X_i$.  
Cover \cite{Cover75} posed some natural questions which 
have been addressed by Bailey \cite{Bailey76}, Ryabko \cite{Ryabko88},
and Ornstein \cite{Ornstein74}, and more recently by
Algoet \cite{Algoet92}, Morvai, Yakowitz, Gy\"{o}rfi \cite{MYGy96}
and Morvai, Yakowitz, Algoet \cite{MYA}. 
Nobel \cite{Nobel97a} has shown that no regression procedure is consistent 
for every bivariate ergodic process, even if one assumes that 
$X_i$ is bounded and $Y_i$ is binary valued.   
A similar negative result for individual sequences is 
established in Theorem \ref{thm2} below.

In order to study regression estimation in a deterministic setting
one must first specify how an individual
sequence $(x_1,y_1), (x_2, y_2), \ldots$ 
can contain information
about a regression function.
In the present paper, following \cite{NobMorKul96}, it is required that  
suitable averages over the sequence are convergent or `stable'.
The deterministic setting of this paper is also in line with other
recent work on individual sequences in information theory, statistics,
and learning theory (cf. \cite{Ziv78,MerFedGut92,HauKivWar94}).
The principal contribution of the paper is 
to show how one may extract asymptotic information from the sequence
in the absence of probabilistic inequalities, mixing conditions, 
rates of convergence, and so on.  The deterministic setting is
described in Section 2 and the principal results of the paper 
are stated in Section 3.
Proofs of the principal results are given in Sections~4 and~5.

\section{The Deterministic Setting}

Let $\mu$ be a probability distribution on $(\real, {\cal B})$, 
and let $m: \real \to \real$ be a function satisfying
$\int |m(x)| \mu(dx) < \infty$.  
Let $\bx = (x_1, x_2, \ldots)$ and $\by = (y_1, y_2, \ldots)$ be 
infinite sequences of real numbers.
For each interval $A \subseteq \real$ define the signed measure
$$
\nu(A)=\int_A m(x) \mu(dx) \, .
$$
For each $n \geq 1$ define the relative frequency
$$
\hmu_n(A)=\frac{1}{n} \sum_{i=1}^n I\{x_i \in A\} \, ,
$$
and the joint sample average 
$$
\hnu_n(A)=\frac{1}{n} \sum_{i=1}^n y_i I\{x_i \in A\} \, .
$$
The sequence $\bx$ will
be said to have {\em limiting distribution} $\mu(\cdot)$ if 
\be
\label{condprob}
\hmu_n(-\infty, t] \to \mu(-\infty,t] \ \mbox{ and } \  
\hmu_n(\{t\}) \to \mu(\{t\}) \   
\mbox{ for every $t \in \real$,}
\ee 
and the pair $(\bx,\by)$ will be said to have 
{\em limiting regression} $m(\cdot)$ if 
\be
\label{condreg}
\hnu_n(-\infty, t] \to \nu(-\infty,t] \ \mbox{ and } \
\hnu_n(\{t\}) \to \nu(\{t\}) \ 
\mbox{ for every $t \in \real$.}
\ee
(Note that the second condition is superfluous in each case
if $\mu$ is non-atomic.)
By minor modification of a standard proof of the Glivenko Cantelli 
Theorem (such as that in 
Pollard \cite{Pol84}), one may show that 
if $\bx$ has limiting distribution $\mu(\cdot)$ then in fact
\be
\label{glivcondprob}
\sup_{A \in \A} |\hmu_n(A) -\mu(A)| \to 0 \, ,
\ee 
where $\A$ is the collection of all intervals of the form $(a,b]$ and 
$(-\infty,b]$ with $a, b \in \real$.

An individual sequence $(\bx,\by)$ satisfying (\ref{condprob}) 
and (\ref{condreg})
will be called {\em stable}. Let $\Omega(\mu,m)$ denote
the set of stable sequences with limiting distribution $\mu$ and
limiting regression $m$.  
Stability concerns only the
asymptotic behavior of $\hmu_n$ and $\hnu_n$, which need not converge
to their respective limits at any particular rate.
No constraints are place on the mechanism by which 
the individual sequences $(\bx,\by)$ are produced.
Note in particular that 
membership of $(\bx,\by)$ in $\Omega(\mu,m)$ is unaffected
if one adds to $\bx$ and $\by$ finite
prefixes $x_1', \ldots, x_k'$ and 
$y_1', \ldots, y_k'$ having the same length.
The next proposition, showing that
the sample paths of ergodic processes
are stable with probability one,
follows easily from Birkhoff's ergodic theorem.

\begin{propref}
\label{birk1}
Let $(X_1,Y_1), (X_2,Y_2), \ldots$ be stationary such that $E|Y| < \infty$.
Then $\{(X_i,Y_i)\}$ is stable with probability one. 
If, in addition, $\{(X_i,Y_i)\}$ is ergodic then 
$({\bf X,Y}) \in \Omega(\mu,m)$
with probability one, 
where $\mu(A)=P(X\in A)$, $m(x)=E(Y|X=x)$, ${\bf X} = (X_1, X_2, \ldots)$ and ${\bf Y} = (Y_1, Y_2, \ldots)$.
\end{propref}

\noindent
{\bf Proof:}
Let $\E$ denote the invariant $\sigma$-algebra. By
Birkoff's pointwise ergodic theorem (cf. Stout~\cite{stout}
Theorem~3.5.6 p. 176), for arbitrary Borel-measurable set $A\subset
\real$ with probability one,  
$$
\hmu_n(A) \to P(X_1\in A|\E)=:\mu_{\E}(A)
$$
and 
$$
\hnu_n(A) \to E(Y_1 I_{\{ X_1\in A\}}|\E)=:\nu_{\E}(A).
$$
If, in addition, $\{(X_i,Y_i)\}$ is ergodic then  $\E$ is the trivial
$\sigma$-algebra and so 
 $$
\hmu_n(A) \to P(X_1\in A)
$$ 
and 
$$
\hnu_n(A)\to E(Y_1 I_{\{ X_1\in A\}}). 
$$
The rest follows from the standard proof of the Glivenko Cantelli 
Theorem (cf.   Pollard \cite{Pol84}). 
$\Box$

\noindent
{\bf Remark 1. }
Note that for individual sequences,  
$$
\hmu_n(-\infty,t] \to \mu(-\infty,t] \ \ \mbox{for all $t\in \real$}
$$ 
does not necessarily imply  
$$
\hmu_n(\{t\}) \to \mu(\{t\}) \ \ \mbox{for all $t\in \real$.}
$$
Indeed, with $\bx = ( {-1 \over 2}, {-1 \over 3}, \dots)$,
$\hmu_n(-\infty,t] = 1$ for $t \geq 0$, while 
$\hmu_n(-\infty,t] \to 0$ for $t < 0$.
Thus the limiting distribution $\mu$ 
should concentrate on the atom $\{0\}$,
but $\hmu_n(\{0\})=0$ for all $n$.

\section{Statement of Principal Results}

Recall that the total variation of a real-valued
function $h$ defined on an interval $(a,b]$
is given by
\[
V(h:a,b) = \sup \sum_{i=1}^n |h(t_i) - h(t_{i-1})| \, ,
\]
where the supremum is taken over all finite ordered sequences
\[
a < t_0 < t_1 < \cdots < t_{n-1} < t_n = b \, .
\]
Let $\pint$ denote the positive intergers.
For each non-decreasing function 
$\alpha: \pint \rightarrow (0,\infty)$, let $\F(\alpha)$ denote
the set of bounded measurable functions $m: \real \to \real$ 
such that $V(m:-i,i) < \alpha(i)$ for all $i \geq 1$.  
Let $\pi_0 = \{\real\}$, and for each $k \geq 1$ let 
$\pi_k$ be the partition of $\real$
consisting of the dyadic intervals 
\[
A_{k,j} = \left( \frac{(j-1)}{2^k} , \frac{j}{2^k} \right]
\ \ -\infty < j < \infty \, .
\]
Let $\pi_k[x]$ denote the unique cell of $\pi_k$ containing $x \in \real$.
Note that $\pi_{k+1}$ refines $\pi_k$, and that for each $x$, 
\[  
\lim_{k \to \infty} \len(\pi_k[x]) = 0 \, , 
\] 
where $\len(A)$ 
denotes the length of an interval $A$.

Let $m\in \F(\alpha)$ be arbitrary. Let $\mu$ denote an arbitrary
probability distribution on $\real$. 
Fix two numerical sequences $\bx$ and $\by$ such that $(\bx,\by)\in
\Omega(\mu,m)$.
For each $k \geq 1$ we define a histogram regression estimate
based on $\pi_k$ and adaptively chosen initial sequences
of $\bx$ and $\by$.
For each $n\ge 1$, $k \geq 0$ define 
\[
\mh_{k,n}(x) = {\hnu_n(\pi_k[x])\over \hmu_n(\pi_k[x])} \, ,
\]
where by convention $0/0 = 0$.  Note that $\mh_{k,n}$ is piecewise
constant on the cells of $\pi_k$.
Let $\tau_0 = 1$ and for each $k \geq 1$ define
\[
\label{tauk}
\tau_k = \min \left\{ n > \tau_{k-1} \, : \, 
           V(\mh_{k,n}: -i,i) < 4 \alpha(i) 
           \ \mbox{for all $1 \leq i \leq k$} \right\}.
\]
By Lemma~\ref{taufinitelemma}, $\tau_k$ is well defined and finite. 
Note that $\tau_k\to\infty$. 
Define the estimate
\[
\mh_k = \mh_{k, \tau_k} \, .
\]
Note that  $\mh_k$ depends only on the pairs
$(x_1,y_1), \ldots, (x_{\tau_k},y_{\tau_k})$.
To create a fixed sample size version of the estimate for $n\ge 1$ let
\[
\kappa_n = \max\{ k \geq 0 : \tau_k \leq n \} 
\] 
and define
\[
\tm_n = \mh_{\kappa_n} .
\] 
The $L_2(\mu)$-consistency of the estimates is established in the 
following theorem.

\vskip.2in

\begin{thm}
\label{thm1} 
Let $\alpha: \pint \rightarrow (0,\infty)$ be a known, non-decreasing
function.  For every $m(\cdot) \in \F(\alpha)$, every
probability distribution $\mu$ 
on $\real$, and every stable pair $(\bx,\by) \in \Omega(\mu,m)$ 
such that the components of $\by$ are bounded, 
\[
\int (\mh_k(x) - m(x))^2 \mu(dx) \to 0 
\ \ \mbox{ and } \ \ 
\int (\tm_n(x) - m(x))^2 \mu(dx) \to 0 \, .
\]
In other words, the estimates
$\tm_n$ and $\mh_k$ are $L_2(\mu)$-consistent.  
\end{thm}

\noindent
{\bf Remark 2.} Definition of $\mh_k$ is based solely on 
$\alpha(\cdot)$ and the given numerical sequences.
In advance of the data, one need only know a bound on the
variation of its limiting regression on the intervals $(-i,i]$.
The limiting distribution $\mu$, the pre-asymptotic
behavior of the individual sequences, and the bound on the $y_i$
need not be known in advance.

\noindent
{\bf Remark 3.} Let $\B(M)$ denote the class of monotone bounded functions $m: \real \to\real $ such that $|m(x)|<M$ for all $x\in \real$. Since $B(M)\subset \F(\alpha)$ with 
$\alpha(n)=2M$ Theorem~\ref{thm1} is applicable to $B(M)$. 
  
\noindent
{\bf Remark 4.} Let $\L(C)$ denote the class of Lipschitz continuous 
functions $m: \real \to \real$ such that 
$|m(x)-m(z)|<C|z-x|$ for all $x,z\in\real$. 
Since $L(M)\subset \F(\alpha)$ with 
$\alpha(n)=2Cn+\epsilon$ where $0<\epsilon<\infty$ is arbitrary,  Theorem~\ref{thm1} is applicable to $\L(C)$. 

\medskip
\noindent
Theorem~\ref{thm1} and Poposition~\ref{birk1} imply the next corollary. 
\medskip

\begin{cor}
\label{cor1}
Let $\alpha: \pint\rightarrow (0,\infty)$ be a known, 
non-decreasing function.
For every stationary ergodic process 
$(X_1,Y_1), (X_2,Y_2), \ldots \in \real \times \real$ such that 
$X_i$ has distribution $\mu$,
$Y$ is bounded with probability one, and
$m(x)=E(Y_i|X_i=x) \in \F(\alpha)$,
\[
\int ( \mh_k(x) - m(x) )^2 \mu(dx) \to 0  \ \ \mbox{and} \ \  
\int (\tm_n(x) - m(x))^2 \mu(dx) \to 0 \, 
\]
with probability one.
\end{cor}

\medskip
\noindent
Theorem~\ref{thm1} and Proposition~\ref{birk1} imply even more. We apply the same notations as in the proof of Proposition~\ref{birk1}.  
\medskip

\begin{cor}
\label{cor2}
Let $\alpha: \pint\rightarrow (0,\infty)$ be a known, 
non-decreasing function. 
\\ Let  
$(X_1,Y_1), (X_2,Y_2), \ldots \in \real \times \real$ be a stationary process such that $Y$ is bounded with probability one.  
Let 
 $m_{\E}:={d\nu_{\E}\over d\mu_{\E}}$, that is, $\nu_{\E}(A)=\int_A m_{\E}(x) \mu_{\E}(dx)$.
Assume that 
$m_{\E}(\cdot) \in \F(\alpha)$ with probability one. Then   
\[
\int ( \mh_k(x) - m_{\E}(x) )^2 \mu_{\E}(dx) \to 0  \ \ \mbox{and} \ \  
\int (\tm_n(x) - m_{\E}(x))^2 \mu_{\E}(dx) \to 0 \, 
\]
with probability one.
\end{cor}

\vskip.2in

\noindent
The conditions in Theorem ~\ref{thm1} cannot be significantly weakened. 
\vskip.2in

\begin{thm}
\label{thm2} Let $\lambda$  denote the uniform distribution on $[0,1]$. 
There is no $L_2(\lambda)$  
consistent regression procedure for the family of
stable sequences $(\bx,\by)$ such that $x_i \in [0,1]$      
has limiting distribution $\lambda$, and $y_i \in \{0,1\}$
has limiting regression $m$ with $V(m:0,1) < \infty$.
\end{thm}

\section{Proof of Theorem \ref{thm1} }

\begin{lem}
\label{taufinitelemma}
Let $\alpha: \pint \rightarrow (0,\infty)$ be a known, non-decreasing
function.  For every $m(\cdot) \in \F(\alpha)$, every
probability distribution
$\mu$ on $\real$, every stable pair $(\bx,\by) \in \Omega(\mu,m)$, and 
for all $k\ge 0$, $\tau_k$ is well defined and finite.  
\end{lem}

\noindent
{\bf Proof:} 
By definition $\tau_0=1$. Hence we may assume $k\ge 1$.  
Let $f$ be any function with bounded variation  $V(f: -i,i)  < \infty$ on 
$(-i,i]$. Define 
\[
(f \circ \pi_k)(x) = \frac{1}{\mu(\pi_k[x])} \int_{\pi_k[x]} f(z) \mu(dz).
\]                   
Note that $f \circ \pi$ is piecewise constant on the cells of $\pi$.

For $f$ non-decreasing it is immediate that $V(f \circ \pi_k :
-i,i)\le V(f : -i,i)$.  If $f$ is not necessarily non-decreasing then
$f(x) = u(x) - v(x)$ where $u(\cdot)$ and $v(\cdot)$ are
non-decreasing, $V(u: -i,i) \leq V(f: -i,i)$ and $V(v: -i,i) \leq 2
V(f: -i,i)$ (cf. Kolmogorov and Fomin \cite{KolFom70}).  It follows
from the definition that $f \circ \pi_k = u \circ \pi_k - v \circ
\pi_k$, and since $u$ and $v$ are non-decreasing, so are $u \circ
\pi_k$ and $v \circ \pi_k$. Therefore
\begin{eqnarray*}
V(f \circ \pi_k : -i, i)&=& V(u \circ \pi_k -v \circ \pi_k: -i, i)\\
&\le&V(u \circ \pi_k : -i, i)+V(v \circ \pi_k : -i, i)\\
&\le& V(u  : -i, i)+V(v  : -i, i)\\
 &\le& 3V(f: -i,i) 
\end{eqnarray*}
as the variation of the sum is less
than the sum of the variations.  
Now note that since $V(m: -i,i)<\alpha(i)$ hence as $n \to \infty$
\begin{eqnarray*}
V(\mh_{k,n} : -i,i) & = &
\sum_{j = -i2^k+1}^{i2^k-1} 
| {\hnu_n(A_{k,j})\over \hmu_n(A_{k,j})} -      
{\hnu_n(A_{k,j+1})\over \hmu_n(A_{k,j+1})}| \\
& \to & 
\sum_{j = -i2^k+1}^{i2^k-1} 
|{  \nu(A_{k,j})\over \mu(A_{k,j})} - 
{\nu(A_{k,j+1})\over \mu(A_{k,j+1})}| \\
& = & V(m \circ \pi_k: -i,i)\\
&\le& 3V(m : -i,i)<4\alpha(i) .
\end{eqnarray*}
Thus $\tau_k$ is well defined and finite. 
$\Box$

\vskip.2in

{\bf Proof of Theorem \ref{thm1}:} 
Fix a sequence $(\bx, \by)$ satisfying the conditions of the theorem.
For each $k \geq 1$ define $g_k(x) = \mh_k(x) - m(x)$.  
It follows from the definition of $\tau_k$ and the assumption that
$m(\cdot) \in \F(\alpha)$ that 
\[
V(g_k:-i,i) \leq V(\mh_k : -i,i) + V(m : -i,i) < 5 \alpha(i) 
\ \mbox{ for $1 \leq i \leq k$.}
\]
Let $D/2 > 1$ be a common bound for $m(\cdot)$ and the elements of
$\by$, so that $|g_k(x)| < D$ for each $x$.

Let $U = \{ u_1, u_2, \dots \}$ be those numbers $u$ 
for which $\mu(\{u\})>0$.  Then $U$ is either finite or 
countably infinite. 
Note that $\mu$ may be decomposed as a sum 
$\mu_d + \mu_c$, where $\mu_d$ is a purely atomic measure supported on
$U$, and $\mu_c$ is non-atomic.
Fix $\epsilon \in (0,1)$.  Let $T \geq 1$ be an integer such that \be
\label{tailofmissmall}
\mu(\{ x : |x| \geq T \}) < \frac{\epsilon}{D^2}
\ee
and let $J \geq 1$ be so large that 
\be
\label{discretetailissmall}
\sum_{i=J+1}^{|U|} \mu(\{u_i\})< {\epsilon\over D^2},
\ee
where $|U|$ denotes the cardinality of $U$.  For $k\ge 1$ define 
$$
\Delta(k) = 
\min\{\mu_c(A) :  A \in \pi_k, \ A \subseteq (-T,T], \ \mu_c(A) > 0 \} 
$$
and 
$$ 
\Theta(k)=\max\{ \mu_c(A) : A \in \pi_k, A \subseteq (-T,T]\}. 
$$
Note that $\Theta(k)\ge \Delta(k) > 0$ for each $k$ and 
that $\Theta(k)$ is a non-increasing 
function of $k$.  Let 
\[
\Theta^* = \lim_{k \to \infty} \Theta(k) \, .
\]

Suppose that $\Theta^* > 0$.  Then there is a sequence of intervals 
$A_k \in \pi_k$ such that $\mu_c(A_k) \geq \Theta^*$ and
$\clos(A_{k+1}) \subseteq \clos(A_k)$ for each $k \geq 1$, where 
$\clos(A)$ denotes the closure of $A$.
As $\len(A_k) \to 0$, $\cap_k \clos(A_k)$ is a singleton $\{x_0\}$.
Continuity of $\mu_c$ implies that 
$\mu_c(\{x_0\}) \geq \Theta^* > 0$, which contradicts the fact that 
$\mu_c$ is non-atomic.  Therefore $\Theta^* = 0$. 
Let $K\ge 1$ be so large that 
\be
\label{intervalissmall}
\Theta(K) < \frac{\epsilon^2}{10 \alpha(T) D^2}.
\ee 
Fix an atom $u \in U$. If $r \leq k$ then
\[
{ \hnu_{\tau_k}(\{u\})+D\hmu_{\tau_k}(\{u\})\over \hmu_{\tau_k}(\pi_r(u))
}
\le 
{ \hnu_{\tau_k}(\pi_k(u))+D\hmu_{\tau_k}(\pi_k(u))\over
\hmu_{\tau_k}(\pi_k(u)) }
\le 
{ \hnu_{\tau_k}(\pi_r(u))+D\hmu_{\tau_k}(\pi_r(u))\over
\hmu_{\tau_k}(\{u\}) }.
\]
As $k$ tends to infinity, stability implies that
$$
\hmu_{\tau_k}(\pi_r(u))\to\mu(\pi_r(u)), \ \ 
\hnu_{\tau_k}(\pi_r(u))\to\nu(\pi_r(u)), \ \ 
\hmu_{\tau_k}(\{u\})\to\mu(\{u\}).  
$$ 
As $r$ tends to infinity, continuity of the measures $\mu$ and $\nu$
implies that 
$$
\mu(\pi_r(u)) \to \mu(\{u\}), \ \ 
\nu(\pi_r(u)) \to \nu(\{u\}). 
$$
From these relations we conclude that 
\be
\label{gkutozero}
\lim_{k \to \infty}  
{ \hnu_{\tau_k}(\pi_k(u))\over \hmu_{\tau_k}(\pi_k(u)) }=
{\nu(\{u\})\over \mu(\{u\}) }.
\ee
By (\ref{glivcondprob}), (\ref{gkutozero}) and (\ref{condreg}) there exists
$K' \geq \max(K,T)$ such that for all indices $k \geq K'$,
\be 
\label{probsareclose}
\sup_{A \in \A} |\hmu_{\tau_k}(A) - \mu(A)| < \frac{\epsilon}{4D}
\Delta(K),
\ee
\be
\label{gkuisarecloseuptoj}
|g_k(u_i)|^2< {\epsilon\over J} \ \ \mbox{for $i=1,\dots,J$},
\ee
and 
\be
\label{estimateisgoodona}
|\int_A \mh_k d\hmu_{\tau_k} - \int_A m d\mu|
\, < \, {\epsilon \over 4} \Delta(K) 
\ee
for every cell $A \in \pi_K$ with $A \subseteq (-T,T]$.

Fix $k \geq K'$, and let $A \in \pi_K$ be such that $\mu(A) > 0$
and $A \subseteq (-T,T]$.  Inequalities
(\ref{probsareclose}) and (\ref{estimateisgoodona}) imply that 
\begin{eqnarray*}
| \int_A g_k(x) \mu(dx) | 
& \leq & 
| \int_A \mh_k d\mu - \int_A \mh_k d\hmu_{\tau_k} |
\, + \, | \int_A \mh_k d\hmu_{\tau_k} - \int_A m d\mu |  \\
& \leq & 
D \sup_{A' \in \A} |\hmu_{\tau_k}(A') - \mu(A')| 
\, + \, \frac{\epsilon}{4} \Delta(K)  \\ 
& \leq & \frac{\epsilon}{2} \Delta(K) \, ,
\end{eqnarray*}
and therefore
\be
\label{intgksmall}
\left| \frac{\int_A g_k(x) \mu(dx)}{\mu(A)} \right|
\leq {\epsilon\over 2}.
\ee
Consider those points
\[
H_k = \{ x \in \real : |g_k(x)| > \epsilon \}
\]
for which $g_k$ exceeds $\epsilon$, and define
\[
\H_k = 
\{ A \in \pi_K \, : \, A \cap H_k \neq \emptyset, 
\, A \subseteq (-T,T], \, 
\mu(A) > 0 \}.
\]
If $A \in \H_k$ then there exists $x \in A$ such that 
$|g_k(x)| > \epsilon$. 
Assume without loss of generality that $g_k(x) > \epsilon$. 
By virtue of (\ref{intgksmall}) there exists $z \in A$ such that 
$g_k(z) \leq \epsilon/2$, and therefore, 
$|g_k(x) - g_k(z)| > \epsilon / 2$ for some $x, z \in A$.
Consequently
\[
{\epsilon \over 2} |\H_k| \leq V(g_k: -T,T) < 5 \alpha(T)
\]
from which follows that 
\be
\label{numberofcalH}
 |\H_k| < {10 \alpha(T) \over \epsilon}.
\ee

Consider now the $L_2(\mu)$ error of $\mh_k$.  From the definition of 
$\H_k$ and inequalities (\ref{numberofcalH}),
(\ref{intervalissmall}), (\ref{discretetailissmall}),  
(\ref{gkuisarecloseuptoj}), and 
(\ref{tailofmissmall})
it follows that
\begin{eqnarray*}
\int |g_k(x)|^2 \mu(dx)
& \leq & \sum_{A \in \H_k} \int_{A} D^2 d\mu_c
\, + \, \sum_{A \in \H_k} \int_{A} |g_k(x)|^2  d\mu_d(x)\\
& + & \sum_{A \notin \H_k, A\subseteq (-T,T]} \int_A \epsilon^2 d\mu
\, + \, \int_{|x| \geq T} D^2 d\mu \\
& \leq &
\epsilon + 
\sum_{i=1}^{J} |g_k(u_i)|^2  +
\sum_{i=J+1}^{|U|} D^2 \mu_d(\{u_i\}) + \epsilon^2  + \epsilon \\
&\leq& 4 \epsilon + \epsilon^2.  
\end{eqnarray*}
Letting $k \to \infty$ and $\epsilon \to 0$ shows that 
$\int |g_k(x)|^2 \mu(dx) \to 0$.
Since $\kappa_n \nearrow \infty$, the $L_2(\mu)$ convergence of
$\tm_n$ to $m$ is immediate. $\Box$

\section{Proof of Theorem \ref{thm2}} 

\noindent
{\bf Proof of Theorem \ref{thm2}:}
For $k \geq 1$ define the $k$'th Rademacher function as 
\[
h_k(x) = 
\left\{ \begin{array}{ll} 
        1 & \mbox{ if $2j 2^{-k} \leq x < (2j+1) 2^{-k}$ for some
                   $0 \leq j < 2^{k-1}$ } \\ 
        0  & \mbox{otherwise ,} \end{array} \right.
\]
and let 
\[
h_0(x)=
\left\{ \begin{array}{ll}
	0.5 & \mbox{ if $x\in [0,1]$}\\
	0 & \mbox{otherwise.} 
\end{array} \right.
\] 
Define $\F_0 = \{ h_0, h_1, h_2, \ldots \}$ and let
$\F_1 = \{ h_1, h_2, \ldots \}$. 
Let $\lambda$ denote the uniform distribution on $[0,1]$. 
We will prove even more than stated in Theorem~\ref{thm2}, namely: 

\vspace{.3cm}
\noindent
{\it There is no $L_2(\lambda)$ consistent regression estimation procedure
for the family
\[ 
\Omega^* = \bigcup_{m \in {\cal F}_0} 
\Omega(\lambda,m) \cap 
\{(\bx,\by): x_n \in [0,1], y_n \in \{0,1\}  
\ \mbox{for all $n \geq 1$} \}.
\]
}
\vspace{.3cm}

This statement says that even for the countable class ${\cal F}_0$ of
regression functions there is no $L_2(\lambda)$ consistent estimation
procedure.  We briefly describe the main idea of the proof.  Let
$\Phi= \{\phi_1, \phi_2, \ldots\}$ be any regression estimation
procedure.  If $\Phi$ fails to be consistent for some sequence
$(\bx,\by) \in \bigcup_{m \in {\cal F}_1} \Omega(\lambda,m)$ with $x_i \in
[0,1]$ and
$y_i \in \{0,1\}$,
there is nothing to prove. 
Assuming then that $\Phi$ is consistent for every such sequence,
we construct a stable sequence $(\bx^*,\by^*)$ 
such that $\phi_n(\cdot : (x_1^*,y_1^*), \ldots,
(x_n^*,y_n^*))$ fails to converge.
The sequence $(\bx^*,\by^*)$ has limiting
distribution $\lambda$ and limiting regression $h_0$.  
It is constructed by `splicing' together longer and longer
blocks of stable sequences 
$(\bx^{(k)},\by^{(k)}) \in \Omega(h_k,\lambda)$.  
When applied to the resulting sequence, the procedure $\Phi$
first produces estimates close to $h_1$; as the sample size is
increased $\Phi$ produces estimates close to $h_2$, then $h_3$,
and so on.  As the $h_i$'s fail to converge, so to do the estimates
$\phi_n(\cdot : (x_1^*,y_1^*), \ldots, (x_n^*,y_n^*))$, $n \geq 1$.

\vspace{1cm}

Note that each $h_j$ is supported on $[0,1]$ and that 
$\int |h_j(x) - h_k(x)|^2 \lambda(dx) = 0.5$ whenever $j \neq k$, and
$j\ge
1$,
$k\ge 1$. 
Let 
\[
\nu_k(A)=\int_{A} h_k(x)\lambda(dx)
\]
and for each finite sequence $(x_1,y_1), \ldots, (x_m,y_m) \in
[0,1]\times \{0,1\}$ let
\[
\Delta(x_1,\dots,x_m) = 
\sup_{A \in {\A}}
\left| \frac{1}{m} \sum_{j=1}^{m} I\{x_j\in A\} - \lambda(A)
\right|=\sup_{A \in {\A}}
\left| \hmu_m(A) - \lambda(A) \right|
\]
and
\begin{eqnarray*}
{\tilde \Delta}_k((x_1,y_1),\dots,(x_m,y_m)) &=& 
\sup_{A \in {\A}}\left| \frac{1}{m} 
\sum_{j=1}^{m} y_j I\{x_j\in A\} - \nu_k(A)\right|\\
&=&\sup_{A \in {\A}}\left| \frac{1}{m} 
\sum_{j=1}^{m}  I\{y_j=1, x_j\in A\} - \nu_k(A)\right|\\
&=&\sup_{A \in {\A}}\left| \hnu_m(A) - \nu_k(A) \right| 
\end{eqnarray*}
where $\A$ is the collection of all intervals of the form $(a,b]$ and
$(-\infty,b]$ with $a, b \in \real$.

A minor modification of a standard proof of the Glivenko Cantelli
Theorem (e.g.\ using the bracketing approach found 
in Pollard \cite{Pol84}) shows that 
\be 
\label{GCT}
\Delta(x_1,\dots,x_m) \to 0 \ \ \mbox{and} \ \   
{\tilde \Delta}_k((x_1,y_1),\dots,(x_m,y_m)) \to 0
\ee
for all  
$({\bf x,y})\in\Omega(\lambda,h_k)\cap \{(\bx,\by): x_n\in (0,1), 
y_n\in \{0,1\} \ \mbox{for all $n\ge 1$}\}$.

Suppose now that $\Phi= \{\phi_1, \phi_2, \ldots\}$ is consistent for
$\F_1$.  For each $k \geq 1$ select a sequence
$$
(\bx^{(k)},\by^{(k)}) = ((x_1^{(k)},y_1^{(k)}), 
(x_2^{(k)}, y_2^{(k)}), \ldots)
$$
such that 
$$
(\bx^{(k)},\by^{(k)}) \in \Omega(h_k,\lambda)\cap \{(\bx,\by): x_n\in
(0,1), 
y_n\in \{0,1\} \ \mbox{for all $n\ge 1$}\}
$$
and 
\be
\label{norepeatition}
x_i^{(k)}=x_j^{(l)} \ \mbox{ if and only if $i=j$, $k=l$}
\ee
(e.g. typical sample sequences from  independent i.i.d.  
time series 
$$
(X^{(k)}_1, h_k(X^{(k)}_1)), (X^{(k)}_2, h_k(X^{(k)}_2)),\dots
$$
where $X^{(k)}_i$ has distribution $\lambda$  
cf. Proposition~\ref{birk1}).  Define  
\be
\label{deflk}
l_k = \min\left\{ L : 
\sup_{m \geq L} 
\Delta(x_1^{(k)}, \ldots, x_m^{(k)})
\leq \frac{1}{k+1}  \right\} 
\ee
\be
\label{deftildelk}
{\tilde l}_k = \min\left\{ L : 
\sup_{m \geq L} 
{\tilde \Delta}_k((x_1^{(k)},y_1^{(k)}), \ldots, (x_m^{(k)},y_m^{(k)}))
\leq \frac{1}{k+1}  \right\}. 
\ee
By (\ref{GCT}), both 
$l_k$ and ${\tilde l}_k$  are finite.
Consider the infinite sequence $(\bx^{(1)},\by^{(1)})$.  
As $h_1 \in \F_1$, and $\Phi$ is consistent for ${\cal F}_1$
by assumption, 
\[
\lim_{n\to\infty}\int | \phi_n(x : (x_1^{(1)},y_1^{(1)}) \ldots,
(x_n^{(1)},y_n^{(1)})) - h_1(x) |^2 \lambda(dx) = 0.
\]
Therefore there is an integer
$n_1 \geq \max(l_2,{\tilde l_2})$ and a corresponding initial segment
$(\bv^{(1)},\bw^{(1)}) = ((x_1^{(1)}, y_1^{(1)}) \ldots, (x_{n_1}^{(1)},y_{n_1}^{(1)}))$ of $(\bx^{(1)},\by^{(1)})$ 
such that 
\[
\int | \phi_{n_1}(x : (\bv^{(1)},\bw^{(1)})) - h_1(x) |^2 \lambda(dx) \leq
\frac{1}{40}
\]
and
\[
\Delta(\bv^{(1)}) \leq \frac{1}{2} \, 
\]
and
\[
{\tilde \Delta}_1((\bv^{(1)},\bw^{(1)})) \leq \frac{1}{2}.
\]
Let $n_0=0$ and let $n_1$ be as defined above. 
Now suppose that for all $1\le j\le k$ 
one has constructed  sequences  $(\bv^{(j)},\bw^{(j)})$ of 
finite length $n_{j}$ in such a way that 
\be
\label{concatenation}
(\bv^{j},\bw^{j})=(v^{(j-1)}_1,w^{(j-1)}_1),\dots,(v^{(j-1)}_{n_{j-1}},w^{(j-1)}_{n_{j-1}}),(x^{(j)}_1,y^{(j)}_1),\dots,(x^{(j)}_{n_{j}-n_{j-1}},y^{(j)}_{n_{j}-n_{j-1}}), 
\ee
\be
\label{close}
\int | \phi_{n_j}(x : (\bv^{(j)},\bw^{(j)})) - h_j(x) |^2 \lambda(dx) \leq
{1\over 40} \, ,
\ee
\be
\label{supint}
\Delta(\bv^{(j)}) \leq (j+1)^{-1}
\ee
\be
\label{supintreg}
{\tilde \Delta}_j((\bv^{(j)},\bw^{(j)})) \leq (j+1)^{-1}
\ee
\be
\label{nk}
n_j \geq j \cdot \max(l_{j+1},{\tilde l}_{j+1}) \, .
\ee
As $(\bv^{(k)},\bw^{(k)})$ is finite, the concatenation 
$$
(v^{(k)}_1,w^{(k)}_1),\dots,(v^{(k)}_{n_{k}},w^{(k)}_{n_{k}}),(x^{(k+1)}_1,y^{(k+1)}_1),(x^{(k+1)}_{2},y^{(k+1)}_{2}),\dots
$$ is contained in
$\Omega(h_{k+1},\lambda)$.  It follows from the
consistency of $\Phi$  that for all large enough  $n$ 
$$
(v^{(k)}_1,w^{(k)}_1),\dots,(v^{(k)}_{n_{k}},w^{(k)}_{n_{k}}),(x^{(k+1)}_1,y^{(k+1)}_1),(x^{(k+1)}_{2},y^{(k+1)}_{2}),\dots, (x_{n-n_k}^{(k+1)},y_{n-n_k}^{(k+1)})
$$
 satisfies (\ref{concatenation}), 
(\ref{close}), (\ref{supint}) and (\ref{supintreg}) with $j$ replaced by $k+1$.
Select $n_{k+1} > n_k$ so large that the same is true of (\ref{nk}).

As $(\bv^{(k+1)},\bw^{(k+1)})$ is an extension of $(\bv^{(k)},\bw^{(k)})$, repeating the
above process indefinitely yields an infinite sequence
$(\bx^*,\by^*)$.  By construction, the functions
$\phi_n (\cdot) = \phi(\cdot : (x_1^*,y_1^*), \ldots, (x_n^*,y_n^*))$ do not
converge in $L_2(\lambda)$.  Indeed, it follows from (\ref{close}) 
and from the inequality $a^2\ge d^2/5-b^2-c^2$ whenever  $(a+b+c)^2=d^2$ that 
\begin{eqnarray*}
\int |\phi_{n_k}(x)-\phi_{n_l}(x)|^2 \lambda(dx) &\ge& 
{1\over 5} \int |\phi_{k}(x)-\phi_{l}(x)|^2 \lambda(dx)\\
&-& \int |\phi_{k}(x)-\phi_{n_k}(x)|^2 \lambda(dx)\\
&-& \int |\phi_{n_l}(x)-\phi_{l}(x)|^2 \lambda(dx)\\ 
&\ge& {1\over 10} -{1\over 40}-{1\over 40}\\
&\ge& {1\over 20}
\end{eqnarray*}
whenever $k \neq l$, $k\ge 1$, $l\ge 1$.

It remains to show that the limiting distribution  of $\bx^*$  
is $\lambda$  and the limiting regression of $(\bx^*,\by^*)$ is $h_0$.  To
this end, 
fix $k > 1$ and let $A \subseteq [0,1]$ 
be an arbitrary interval. 
It is easily verified that 
\be
\label{nuknu0close}
|\nu_k(A) - \nu_0(A)| \leq 2^{-k+1} \le \frac{2}{k}. 
\ee
Let $\hat{\mu}_n(A)$ and $\hat{\nu}_n(A)$ be 
 evaluated on  $((x_1^*,y_1^*), \ldots, (x_n^*,y_n^*))$, and for each 
$1 \leq r \leq n_{k+1} - n_k$ define
\[
\hat{\nu}_{r,k}'(A) = \frac{1}{r} \sum_{j=n_k+1}^{n_k+r} y_i^* I\{x_i^*\in A\}.
\]
The equation 
\[
\hat{\nu}_{n_k+r}(A) = 
\frac{n_k}{n_k+r} \cdot \hat{\nu}_{n_k}(A) + 
\frac{r}{n_k+r} \cdot \hat{\nu}_{r,k}'(A) 
\]
implies the bound
\begin{eqnarray*}
\label{ratreg}
|\hat{\nu}_{n_k+r}(A) - \nu_0(A)| & \leq &
\frac{n_k}{n_k+r} \cdot |\hat{\nu}_{n_k}(A) - \nu_0(A)| \ + \ 
\frac{r}{n_k+r} \cdot |\hat{\nu}_{r,k}'(A) - \nu_0(A)| \\
& \deq & I + II.
\end{eqnarray*}
By virtue of 
(\ref{supintreg}) and (\ref{nuknu0close})
\[
I \, \leq \,
|\hat{\nu}_{n_k}(A) - \nu_{k}(A)| + |\nu_0(A) - \nu_{k}(A)| 
\, \leq \, \frac{1}{k+1} + \frac{2}{k}. 
\]
If $n_{k+1}-n_k \geq r \geq {\tilde l}_{k+1}$ then 
by 
(\ref{deftildelk})
\begin{eqnarray*}
{\tilde \Delta}_{k+1}((x_{n_k+1}^*,y_{n_k+1}^*), \ldots, 
(x_{n_k+r}^*,y_{n_k+r}^*))
&=& {\tilde \Delta}_{k+1}((x_1^{(k+1)},y_1^{(k+1)}), 
\ldots, (x_r^{(k+1)},y_r^{(k+1)}))\\
&\leq&  \frac{1}{k+2}
\end{eqnarray*}
and therefore
\[
II \, \leq \,
|\hat{\nu}_{r,k}'(A) - \nu_{(k+1)}(A)| \, + \,
|\nu_{(k+1)}(A) - \nu_0(A)| \, \leq \,
\frac{1}{k+2} + \frac{2}{k+1} \, .
\]
On the other hand, if $0<r < {\tilde l}_{k+1}$ then (\ref{nk}) implies
that
\[
II\, \leq \, \frac{2r}{n_k + r} 
\, \leq \, \frac{2r}{kr + r}
\, = \, \frac{2}{k+1} \, .
\]
These bounds ensure that, since $A$ was an arbitrary interval, 
 \[
\max\left\{ \sup_{A\in\A}|\hat{\nu}_n(A) - \nu_0(A)| : 
n_k < n \leq n_{k+1} \right\}
\leq \frac{6}{k}
\]
and consequently,
\[
\lim_{n \to \infty} \sup_{A\in\A}|\hat{\nu}_n(A) - \nu_0(A)| = 0. 
\]
A similar (in fact, easier) analysis  establishes  
\[
\lim_{n \to \infty} \sup_{A\in\A}|\hat{\mu}_n(A) - \lambda(A)| = 0. 
\]
Finally, by~(\ref{norepeatition}), for all $t\in\real$ 
$$
\hat{\mu}_n(\{t\}) \to \lambda(\{t\}) = 0 \ \ \mbox{and} \ \ 
\hat{\nu}_n(\{t\}) \to \nu_0(\{t\}) = 0. \ \ \Box
$$

\vskip.3in

\noindent
{\Large\bf Acknowledgements}

\noindent 
The first author wishes to thank Andr\'{a}s Antos for helpful discussions.



\small{

}

\end{document}